\documentclass[10pt]{article}
\usepackage{latexsym}

\newtheorem{lemma}{Lemma}
\newtheorem{prop}{Proposition}

\newtheorem{theorem}{Theorem}
\newtheorem{coro}{Corollary}

\def\V{{\cal V}}
\def\A{{\cal A}}
\def\B{{\cal B}}

\def\L{L}
\def\He{H_1}
\def\Hz{H_2}
\def\Hd{H_3}
\def\deg{{\rm deg}}
\def\chib{\overline{\chi}}
\def\D21{{D_\alpha}}
\def\Dvog{{\tilde{D}(2,1)}}
\def\WDroa{W_{\D21,\rho_\alpha,\omega_\alpha}}
\def\id{\rm id}

\def\End{\rm End}
\def\diag{\rm diag}

\def\mathbb{\bf}

\newcommand{\Z}{{\mathbb Z}}

\newcommand{\C}{{\mathbb C}}

\def\DottedCircle{
\bezier{4}(0.966,-0.259)(1.04,0)(0.966,0.259)
\bezier{4}(0.966,0.259)(0.897,0.518)(0.707,0.707)
\bezier{4}(0.707,0.707)(0.518,0.897)(0.259,0.966)
\bezier{4}(0.259,0.966)(0,1.04)(-0.259,0.966)
\bezier{4}(-0.259,0.966)(-0.518,0.897)(-0.707,0.707)
\bezier{4}(-0.707,0.707)(-0.897,0.518)(-0.966,0.259)
\bezier{4}(-0.966,0.259)(-1.04,0)(-0.966,-0.259)
\bezier{4}(-0.966,-0.259)(-0.897,-0.518)(-0.707,-0.707)
\bezier{4}(-0.707,-0.707)(-0.518,-0.897)(-0.259,-0.966)
\bezier{4}(-0.259,-0.966)(0,-1.04)(0.259,-0.966)
\bezier{4}(0.259,-0.966)(0.518,-0.897)(0.707,-0.707)
\bezier{4}(0.707,-0.707)(0.897,-0.518)(0.966,-0.259)
}
\def\Endpoint[#1]{
\ifcase#1
\fi}
\def\Arc[#1]{
\thicklines                     
\ifcase#1
\bezier{25}(0.966,-0.259)(1.04,0)(0.966,0.259)
\or
\bezier{25}(0.966,0.259)(0.897,0.518)(0.707,0.707)
\or
\bezier{25}(0.707,0.707)(0.518,0.897)(0.259,0.966)
\or
\bezier{25}(0.259,0.966)(0,1.04)(-0.259,0.966)
\or
\bezier{25}(-0.259,0.966)(-0.518,0.897)(-0.707,0.707)
\or
\bezier{25}(-0.707,0.707)(-0.897,0.518)(-0.966,0.259)
\or
\bezier{25}(-0.966,0.259)(-1.04,0)(-0.966,-0.259)
\or
\bezier{25}(-0.966,-0.259)(-0.897,-0.518)(-0.707,-0.707)
\or
\bezier{25}(-0.707,-0.707)(-0.518,-0.897)(-0.259,-0.966)
\or
\bezier{25}(-0.259,-0.966)(0,-1.04)(0.259,-0.966)
\or
\bezier{25}(0.259,-0.966)(0.518,-0.897)(0.707,-0.707)
\or
\bezier{25}(0.707,-0.707)(0.897,-0.518)(0.966,-0.259)
\fi}
\def\DottedArc[#1]{
\ifcase#1
\bezier{4}(0.966,-0.259)(1.04,0)(0.966,0.259)
\or
\bezier{4}(0.966,0.259)(0.897,0.518)(0.707,0.707)
\or
\bezier{4}(0.707,0.707)(0.518,0.897)(0.259,0.966)
\or
\bezier{4}(0.259,0.966)(0,1.04)(-0.259,0.966)
\or
\bezier{4}(-0.259,0.966)(-0.518,0.897)(-0.707,0.707)
\or
\bezier{4}(-0.707,0.707)(-0.897,0.518)(-0.966,0.259)
\or
\bezier{4}(-0.966,0.259)(-1.04,0)(-0.966,-0.259)
\or
\bezier{4}(-0.966,-0.259)(-0.897,-0.518)(-0.707,-0.707)
\or
\bezier{4}(-0.707,-0.707)(-0.518,-0.897)(-0.259,-0.966)
\or
\bezier{4}(-0.259,-0.966)(0,-1.04)(0.259,-0.966)
\or
\bezier{4}(0.259,-0.966)(0.518,-0.897)(0.707,-0.707)
\or
\bezier{4}(0.707,-0.707)(0.897,-0.518)(0.966,-0.259)
\fi}
\def\Chord[#1,#2]{
\thinlines
\ifnum#1>#2\Chord[#2,#1]
\else\ifnum#1<#2
\ifcase#1
\ifcase#2
\or\qbezier(1,0)(0.516,0.138)(0.866,0.5)
\or\qbezier(1,0)(0.45,0.26)(0.5,0.866)
\or\qbezier(1,0)(0.327,0.327)(0,1)
\or\qbezier(1,0)(0.179,0.311)(-0.5,0.866)
\or\qbezier(1,0)(0.0536,0.2)(-0.866,0.5)
\or\put(1, 0){\line(-2, 0){2}}
\or\qbezier(1,0)(0.0536,-0.2)(-0.866,-0.5)
\or\qbezier(1,0)(0.179,-0.311)(-0.5,-0.866)
\or\qbezier(1,0)(0.327,-0.327)(0,-1)
\or\qbezier(1,0)(0.45,-0.26)(0.5,-0.866)
\or\qbezier(1,0)(0.516,-0.138)(0.866,-0.5)
\fi
\or\ifcase#2\or
\or\qbezier(0.866,0.5)(0.378,0.378)(0.5,0.866)
\or\qbezier(0.866,0.5)(0.26,0.45)(0,1)
\or\qbezier(0.866,0.5)(0.12,0.446)(-0.5,0.866)
\or\qbezier(0.866,0.5)(0,0.359)(-0.866,0.5)
\or\qbezier(0.866,0.5)(-0.0536,0.2)(-1,0)
\or\put(0.866, 0.5){\line(-5, -3){1.73}}
\or\qbezier(0.866,0.5)(0.146,-0.146)(-0.5,-0.866)
\or\qbezier(0.866,0.5)(0.311,-0.179)(0,-1)
\or\qbezier(0.866,0.5)(0.446,-0.12)(0.5,-0.866)
\or\qbezier(0.866,0.5)(0.52,0)(0.866,-0.5)
\fi
\or\ifcase#2\or\or
\or\qbezier(0.5,0.866)(0.138,0.516)(0,1)
\or\qbezier(0.5,0.866)(0,0.52)(-0.5,0.866)
\or\qbezier(0.5,0.866)(-0.12,0.446)(-0.866,0.5)
\or\qbezier(0.5,0.866)(-0.179,0.311)(-1,0)
\or\qbezier(0.5,0.866)(-0.146,0.146)(-0.866,-0.5)
\or\put(0.5, 0.866){\line(-3, -5){1}}
\or\qbezier(0.5,0.866)(0.2,-0.0536)(0,-1)
\or\qbezier(0.5,0.866)(0.359,0)(0.5,-0.866)
\or\qbezier(0.5,0.866)(0.446,0.12)(0.866,-0.5)
\fi
\or\ifcase#2\or\or\or
\or\qbezier(0,1.)(-0.138,0.516)(-0.5,0.866)
\or\qbezier(0,1.)(-0.26,0.45)(-0.866,0.5)
\or\qbezier(0,1.)(-0.327,0.327)(-1,0)
\or\qbezier(0,1.)(-0.311,0.179)(-0.866,-0.5)
\or\qbezier(0,1.)(-0.2,0.0536)(-0.5,-0.866)
\or\put(0, 1){\line(0, -2){2}}
\or\qbezier(0,1.)(0.2,0.0536)(0.5,-0.866)
\or\qbezier(0,1.)(0.311,0.179)(0.866,-0.5)
\fi
\or\ifcase#2\or\or\or\or
\or\qbezier(-0.5,0.866)(-0.378,0.378)(-0.866,0.5)
\or\qbezier(-0.5,0.866)(-0.45,0.26)(-1,0)
\or\qbezier(-0.5,0.866)(-0.446,0.12)(-0.866,-0.5)
\or\qbezier(-0.5,0.866)(-0.359,0)(-0.5,-0.866)
\or\qbezier(-0.5,0.866)(-0.2,-0.0536)(0,-1)
\or\put(-0.5, 0.866){\line(3, -5){1}}
\or\qbezier(-0.5,0.866)(0.146,0.146)(0.866,-0.5)
\fi
\or\ifcase#2\or\or\or\or\or
\or\qbezier(-0.866,0.5)(-0.516,0.138)(-1,0)
\or\qbezier(-0.866,0.5)(-0.52,0)(-0.866,-0.5)
\or\qbezier(-0.866,0.5)(-0.446,-0.12)(-0.5,-0.866)
\or\qbezier(-0.866,0.5)(-0.311,-0.179)(0,-1)
\or\qbezier(-0.866,0.5)(-0.146,-0.146)(0.5,-0.866)
\or\put(-0.866, 0.5){\line(5, -3){1.73}}
\fi
\or\ifcase#2\or\or\or\or\or\or
\or\qbezier(-1,0)(-0.516,-0.138)(-0.866,-0.5)
\or\qbezier(-1,0)(-0.45,-0.26)(-0.5,-0.866)
\or\qbezier(-1,0)(-0.327,-0.327)(0,-1)
\or\qbezier(-1,0)(-0.179,-0.311)(0.5,-0.866)
\or\qbezier(-1,0)(-0.0536,-0.2)(0.866,-0.5)
\fi
\or\ifcase#2\or\or\or\or\or\or\or
\or\qbezier(-0.866,-0.5)(-0.378,-0.378)(-0.5,-0.866)
\or\qbezier(-0.866,-0.5)(-0.26,-0.45)(0,-1)
\or\qbezier(-0.866,-0.5)(-0.12,-0.446)(0.5,-0.866)
\or\qbezier(-0.866,-0.5)(0,-0.359)(0.866,-0.5)
\fi
\or\ifcase#2\or\or\or\or\or\or\or\or
\or\qbezier(-0.5,-0.866)(-0.138,-0.516)(0,-1)
\or\qbezier(-0.5,-0.866)(0,-0.52)(0.5,-0.866)
\or\qbezier(-0.5,-0.866)(0.12,-0.446)(0.866,-0.5)
\fi
\or\ifcase#2\or\or\or\or\or\or\or\or\or
\or\qbezier(0,-1.)(0.138,-0.516)(0.5,-0.866)
\or\qbezier(0,-1.)(0.26,-0.45)(0.866,-0.5)
\fi
\or\ifcase#2\or\or\or\or\or\or\or\or\or\or
\or\qbezier(0.5,-0.866)(0.378,-0.378)(0.866,-0.5)
\fi\fi\fi\fi}
\def\FullChord[#1,#2]{
\Endpoint[#1]
\Endpoint[#2]
\Arc[#1]
\Arc[#2]
\Chord[#1,#2]
}
\def\EndChord[#1,#2]{
\Endpoint[#1]
\Endpoint[#2]
\Chord[#1,#2]
}
\def\Picture#1{
\begin{picture}(2,1)(-1,-0.167)
#1
\end{picture}
}
\def\DottedChordDiagram[#1,#2]{
\Picture{\DottedCircle \FullChord[#1,#2]}
}
\begin{document}
{\parindent0cm
\begin{center}
\bf\large
ON VASSILIEV INVARIANTS NOT COMING FROM SEMISIMPLE LIE ALGEBRAS
\end{center}

\smallskip

\begin{center}
\small Jens Lieberum

\scriptsize
Institut de Recherche Math\'{e}matique Avanc\'{e}e,
Universit\'{e} Louis Pasteur - C.N.R.S.

7 rue Ren\'{e} Descartes,
67084 Strasbourg Cedex, France

Email: lieberum@math.u-strasbg.fr
\end{center}

\vspace{1in}

\begin{abstract}
We prove a refinement of Vogel's statement that the Vassiliev invariants of
knots
coming from semisimple Lie algebras do not generate all Vassiliev invariants.
This refinement takes into account the second grading on the Vassiliev
invariants induced by cabling of knots.
As an application we get an amelioration of the
actually known lower bounds for the dimensions of the space of
Vassiliev invariants.
\end{abstract}


\section*{Introduction}

Vogel found a
primitive element in the Hopf algebra of chord diagrams that cannot be
detected by semisimple Lie superalgebras with Casimir (\cite{Vog}).
This element is shown to be non-trivial thanks to
an injection of a certain algebra $\Lambda$ into the space of
Chinese characters with two univalent vertices.
As a consequence of this
result not all Vassiliev invariants of knots
come from semisimple Lie superalgebras with Casimir.

In this paper we will use Vogel's methods to construct Chinese
characters that cannot be detected by
semisimple Lie algebras. In contrast to Vogel's element, these Chinese
characters may have any given even number of univalent vertices.
To prove that they are non-trivial, we evaluate them with a weight system
associated to the Lie superalgebra $D(2,1,\alpha)$.
For these computations we will make use of the work of Kricker
(\cite{Kr2}) that simplified and generalized a part of
Bar-Natan and Garoufalidis's
proof of the Melvin-Morton-Rozansky conjecture (\cite{BNG}).
As a consequence of our result for every even~$k$ there is a
Vassiliev invariant that is an eigenvector for the $n$-th cabling of knots
modulo invariants of lower degree
with eigenvalue~$n^k$, but that is not a linear combination of Vassiliev
invariants coming from semisimple Lie algebras.

\section{The Main Result}

Let the ground field always be the field $\C{}$ of complex numbers.
We refer to \cite{BN1} for the definition of the space $\A$ generated by
Chinese character
diagrams and the so-called (STU)-relations, of the space $\B$ generated by
Chinese characters and the (IHX)- and (AS)-relations, and of the isomorphism
$\chib:\B\longrightarrow\A$.
Let us denote the universal weight system associated to a Lie
superalgebra~$\L$ with Casimir element~$\omega$ by
$W_{\L,\omega}:\A\longrightarrow U(\L)$.
It takes values in the center of the universal enveloping algebra $U(L)$
of~$L$.
If $\rho:U(\L)\longrightarrow \End(M)$ is an irreducible finite-dimensional
representation of $\L$, then we get a $\C$-algebra morphism
$W_{\L,\rho,\omega}:\A\longrightarrow \C$
such that $\rho\circ W_{\L,\omega}(.)=W_{\L,\rho,\omega}(.) \id_M$. Important
classes of weight systems consist of linear combinations of the maps
$W_{\L,\rho,\omega}$ where $\L$ is a Lie (super)algebra
with Casimir element $\omega$. Now we can formulate our
main result.

\begin{theorem}\label{ineq}
For $d=15$ or $d\geq 17$ and all even $k\geq 2$
there exists a non-zero linear combination $D_{k+d,k}\in \cal B$ of
connected Chinese characters of degree~$k+d$ and with~$k$ univalent
vertices such that $w(D_{k+d,k})=0$ for all weight
systems~$w$ coming from semisimple Lie algebras.
\end{theorem}

The theorem
implies a conjecture on marked surfaces
stated in Remark~2 of~\cite{BN2}. Therefore it
allows to improve the actually known lower bounds of the space of
Vassiliev invariants~(\cite{BN2}, \cite{Kne}).
The special case~$k=2$ of Theorem~\ref{ineq}
is implied by~\cite{Vog}, Theorem~7.4
(the bounds for~$d$ come from Proposition~\ref{calc} of this
paper). In the case~$k=2$, the
result holds more generally for all semisimple Lie superalgebras with Casimir.
It would be nice
to prove Theorem~\ref{ineq} in this generality also for~$k\geq 4$.

Let us give an equivalent formulation of Theorem~\ref{ineq} on the level of
Vassiliev invariants.
The $n$-th cabling of framed
knots (see \cite{BN1}, Definition 3.13) induces an algebra
endomorphism $c_n$ of the space of Vassiliev invariants of framed
knots~$\V$.
Denote by~$\V_{m,k}$ the space of canonical Vassiliev invariants~$v$
of degree~$m$ such that $c_n(v)-n^k v$ is a Vassiliev invariant of degree
$m-1$. A Vassiliev invariant is called primitive, if it is additive for
the connected sum of framed knots.
In the proof of the following corollary we recall the relation between
$\V_{m,k}$ and two gradings on $\cal B$.

\begin{coro}
For $d=15$ or $d\geq 17$ and all even $k\geq 2$ there exists
a primitive Vassiliev invariant $v\in\V_{k+d,k}$
that is not a linear combination of Vassiliev invariants coming from
semisimple Lie algebras.
\end{coro}
{\bf Proof:}
Let $D_{k+d,k}$ be the element of $\B$ considered in Theorem \ref{ineq}.
Choose some complement $N$ of $\C\cdot D_{k+d,k}$ in $\cal B$ that
respects the two gradings on $\cal B$ given by half of the number of
vertices and the
number of univalent vertices of Chinese characters and such that $\chib(N)$
contains the space of decomposable elements of $\A$. Define a weight system
$w$ that vanishes on $\chib(N)$ and is $1$ on $\chib(D_{k+d,k})$.
Then $w$ is a primitive element of $\A^*$.
So the weight system~$w$ integrates to a primitive Vassiliev invariant~$v$
of degree~$k+d$.
By Theorem~3 of \cite{Kr1} the element $\chib(D_{k+d,k})$ is an eigenvector
for the maps $\psi^n$ from \cite{BN1} with eigenvalue $n^k$ since
$D_{k+d,k}$ has $k$ univalent vertices.
Now Exercise~3.14 of \cite{BN1}
implies $v\in\V_{k+d,k}$.
$\Box$

\smallskip

The rest of the paper is devoted to the proof of Theorem~\ref{ineq}.

\section{Sequences of Weight Systems}

Let $L$ be one of the following simple finite-dimensional Lie
superalgebras (see \cite{Hum},~\cite{Kac}):

$$
{\rm sl}(m,n)\, (m\not=n), {\rm osp}(m,2n), E_6,E_7,E_8,F_4,G_2,
D(2,1,\alpha), F(4), G(3).
$$

Choose a Cartan subalgebra~$H$ of~$L$
and a system of positive roots~$\Delta^+ \subset H^*$.
For a root~$\beta\not=0$ of~$L$ let

\begin{equation}
L_\beta=\left\{v\in L\,\mid\, \mbox{$[h,v]=\beta(h)v$ for all $h\in H$}
\right\}.
\end{equation}

Let $x_\beta$ be a generator of the one-dimensional space $L_\beta$
($\beta\in\Delta^+$).
For $\lambda\in H^*$ let $I_\lambda$
be the left ideal in $U(L)$ generated by the elements $x_\beta$ for
$\beta\in\Delta^+$ and
$h-\lambda(h)$ for $h\in H$.
The module $V_\lambda:=U(L)/I_\lambda$ is called the
Verma module of weight
$\lambda$. Denote the image of $1$ in $V_\lambda$ by $v_\lambda$.
Let $\{y_1,\ldots,y_r\}$ be a basis of
$\bigoplus_{\beta\in \Delta^+} L_{-\beta}$.
Then the Poincar\'{e}-Birkhoff-Witt basis
theorem for $L$ (see \cite{Kac}, Section~3) implies that

\begin{equation}
\left\{y_1^{e_1}\ldots y_r^{e_r}\cdot v_\lambda\ \mid\
\deg\, y_i=\overline{0}\Rightarrow e_i \geq 0,
\quad \deg\, y_i=\overline{1}\Rightarrow e_i \in\{0,1\}\right\}
\end{equation}

is a basis of $V_\lambda$.
The vector $v_\lambda\in V_\lambda$ is a basis of
the one-dimensional space of vectors of weight $\lambda$ in $V_\lambda$;
it generates $V_\lambda$ as a $U(L)$-module.

By Proposition~5.3 c) of \cite{Kac} there exists an invariant,
supersymmetric, regular bilinear form $<,>$ on $L$.
Let $\omega\in L\otimes L$ be the Casimir element associated to $<,>$.
Since for $a\in \A$ the element $W_{L,\omega}(a)$ belongs
to the center of $U(L)$, the endomorphism of~$V_\lambda$ induced by
$W_{L,\omega}(a)$ is a scalar multiple of the identity.

Let $\rho:U(L)\longrightarrow \End(V)$ be a
finite-dimensional irreducible representation with highest weight~$\lambda$.
Then~$V$ is a quotient of~$V_\lambda$ and we may determine
$W_{L,\rho,\omega}(a)$ with computations in~$V_\lambda$.
This enables us to compare weight systems when the weight of the
representation varies because we may use one and the same index set for
the bases of the involved modules.
Lemma~3.9 of \cite{Kr2} implies the following lemma.

\begin{lemma}\label{poly1}
Let $\rho_n$ be a sequence of finite-dimensional irreducible
representations of $L$ with highest weight $n\lambda$.
Then for a Chinese character diagram $D$ the values
$W_{L,\rho_n,\omega}(D)$ depend polynomially on $n$.
If $D$ has $k$ vertices lying on the
oriented circle, then the degree of this polynomial is $\leq k$.
\end{lemma}

The Alexander-Conway weight systems of degree $k$ are those that vanish
on diagrams that have at most $k-1$ vertices on the oriented circle.
By Lemma~\ref{poly1} the coefficient of $n^k$ in
$W_{L,\rho_n,\omega}:\A_k\longrightarrow \C$ is in the algebra of
Alexander-Conway weight systems.
The goal of \cite{Kr2} was to find
an explicit expression for this highest coefficient.

For $v\in L_\beta$
we will denote $\deg\, v\in\{\overline{0},\overline{1}\}$ by $\deg\,\beta$.
Since the restriction of the bilinear form $<,>$ to $H$ is nondegenerate,
we identify $H$ and $H^*$ and get a bilinear form
on $H^*$, also denoted by $<,>$. Then Theorem~3.11 of \cite{Kr2} states the
following.

\begin{theorem}[Kricker]\label{highestcoeff}
For even $k\geq 2$ let $T_k$ be the diagram of degree $k$
from the right-hand side of Figure~\ref{SunDiag}.
Let $\rho_n$ be a sequence of finite-dimensional irreducible representations
of $L$ with highest weight $n\lambda$. Then the coefficient of $n^k$ in
$W_{L,\rho_n,\omega}(T_k)$ is equal to
\nopagebreak
$$
2 \sum_{\beta\in\Delta^+}
(-1)^{\deg\,\beta} <\lambda,\beta>^k.
$$
\end{theorem}

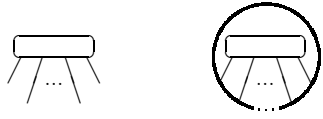
\begin{figure}[!h]
$$
\begin{picture}(2,1.6)(-1,-0.6)
\thinlines
\put(-0.866,-0.5){\line(1,2){0.25}}
\put(-0.5,-0.866){\line(1,3){0.289}}
\put(0.866,-0.5){\line(-1,2){0.25}}
\put(0.5,-0.866){\line(-1,3){0.289}}
\put(-1,-1){\makebox(2,1){$\scriptstyle\ldots$}}
\put(0,0.2){\oval(1.5,0.4)}
\end{picture} \qquad \qquad
\begin{picture}(2,1.6)(-1,-0.6)
\DottedCircle
\Arc[0]\Arc[1]\Arc[2]\Arc[3]\Arc[4]\Arc[5]
\Arc[6]\Arc[7]\Arc[8]\Arc[10]\Arc[11]
\Endpoint[7]\Endpoint[8]\Endpoint[10]\Endpoint[11]
\thinlines
\put(-0.866,-0.5){\line(1,2){0.25}}
\put(-0.5,-0.866){\line(1,3){0.289}}
\put(0.866,-0.5){\line(-1,2){0.25}}
\put(0.5,-0.866){\line(-1,3){0.289}}
\put(-1,-1){\makebox(2,1){$\scriptstyle\ldots$}}
\put(0,0.2){\oval(1.5,0.4)}
\end{picture}
$$
\caption{The diagrams $S_k$ and $T_k$ having $k$ legs}\label{SunDiag}
\end{figure}

\section{Weight Systems Coming from $D(2,1,\alpha)$}

Because of frequent use we will denote the Lie superalgebra $D(2,1,\alpha)$
simply by~$\D21$.
Let us give an explicit description of it.
For $\alpha\in\C\setminus\{0,-1\}$ there exists a basis $(E_i,H_i,F_i)$
($i=1,2,3$) of the even part $(\D21)_{\overline{0}}$ and a basis
$(v_{\epsilon_1 \epsilon_2 \epsilon_3})$
($\epsilon_i\in\{\pm 1\}$) of the odd part $(\D21)_{\overline{1}}$ such that

\begin{itemize}
\item $[H_i,E_j]=2\delta_{ij}E_j, \quad [H_i,F_j]=-2\delta_{ij}F_j, \quad
[H_i,H_j]=0, \quad [E_i,F_j]=\delta_{ij}H_j$,

\item $[H_i,v_{\epsilon_1 \epsilon_2 \epsilon_3}]=
\epsilon_i v_{\epsilon_1 \epsilon_2 \epsilon_3}$,\\
$[E_i,v_{\epsilon_1 \epsilon_2 \epsilon_3}]=
\delta_{-1 \epsilon_i}v_{\gamma_1 \gamma_2 \gamma_3} \quad$ with
$\gamma_\nu=\epsilon_\nu$ for $\nu\not=i$ and $\gamma_i=1$,\\
$[F_i,v_{\epsilon_1 \epsilon_2 \epsilon_3}]=
\delta_{1 \epsilon_i}v_{\gamma_1 \gamma_2 \gamma_3} \quad$ with
$\gamma_\nu=\epsilon_\nu$ for $\nu\not=i$ and $\gamma_i=-1$,

\item $[v_{\epsilon_1 \epsilon_2 \epsilon_3},v_{\gamma_1 \gamma_2 \gamma_3}]=
(\alpha+1) \beta_2\beta_3 G_1(\epsilon_1,\gamma_1)
- \beta_1\beta_3 G_2(\epsilon_2,\gamma_2)
- \alpha \beta_1\beta_2 G_3(\epsilon_3,\gamma_3)$,


where $G_i(1,1)=-E_i$, $G_i(1,-1)=G_i(-1,1)=H_i/2$, $G_i(-1,-1)=F_i$ 
and~$\beta_i=\epsilon_i\delta_{-\epsilon_i\gamma_i}$.
\end{itemize}

In \cite{Vog} a Lie superalgebra $\Dvog$ over the ring $R=\Z[a,b,c]/(a+b+c)$
is considered. Mapping $(a,b,c)$ to $(-\alpha-1,1,\alpha)$ we turn $\C$ into
an $R$-module and recover
$\D21$ as $\Dvog \otimes_R \C$.
Let us choose generators

\begin{equation}
e_1=v_{1 -1 -1},\quad
h_1=\left((\alpha+1)H_1+H_2+\alpha H_3\right)/2,\quad
f_1=v_{-1 1 1}
\end{equation}

and $e_i=E_i$, $h_i=H_i$, $f_i=F_i$ for $i\in\{2,3\}$
of $\D21$. Then one can check that our definition of $\D21$ agrees with
the one given in Section~5 of \cite{Kac}.

We may choose a bilinear form $<,>_\alpha$ on $\D21$ as in
Proposition~5.3 of \cite{Kac}, such that the matrix of the restriction
of $<,>_\alpha$ to the Cartan subalgebra $H$ with basis $(\He,\Hz,\Hd)$ is
$\diag(2(1+\alpha)^{-1},-2,-2\alpha^{-1})$.
Let $\omega_i=E_i\otimes F_i+H_i\otimes H_i/2+F_i\otimes E_i$ and let

\begin{equation}
\pi=\sum_{(\epsilon_i)\in\{\pm 1\}^3} \epsilon_1\epsilon_2\epsilon_3
v_{\epsilon_1 \epsilon_2 \epsilon_3}
\otimes  v_{-\epsilon_1 -\epsilon_2 -\epsilon_3}.
\end{equation}

Then

\begin{equation}
\omega_\alpha=(1+\alpha)\omega_1 -\omega_2-\alpha \omega_3+\pi
\end{equation}

is the
Casimir element in $\D21\otimes\D21$ corresponding to $<,>_\alpha$.
The Casimir element~$\Omega$
from Lemma~6.12 of \cite{Vog} is mapped to $\omega_\alpha$ by the
specialization of the parameters $a,b,c$.

Let $(H_1^*,H_2^*,H_3^*)$ be the basis of $H^*$ dual to $(H_1,H_2,H_3)$.
We choose

\begin{equation}\label{D21deltaplus}
\Delta^+=\left\{2\He^*,2\Hz^*,2 \Hd^*\right\}\cup
\left\{\sum_{i=1}^3 \epsilon_i H_i^*\mid (\epsilon_i)\in
\{(1,\pm 1, \pm 1)\}\right\}
\end{equation}

as a positive root system of $\D21$.

\begin{lemma}\label{poly2}
Let $\rho_{\alpha,n}$ be a sequence of finite-dimensional irreducible
representations of $\D21$ with highest weight $n\lambda$.
Then for a Chinese character diagram $D$ the values
$W_{\D21,\rho_{\alpha,n},\omega_\alpha}(D)$
depend polynomially on $n$ and $\alpha$.
If $D$ has $k$ vertices lying on the
oriented circle, then the degree in $n$ of this polynomial is $\leq k$.
\end{lemma}
{\bf Proof:}
The coefficients in $\omega_\alpha$ and the coefficients of the bracket
of two basis vectors of $\D21$ are polynomials of degree $\leq 1$
in $\alpha$.
Now arguments from the proof of \cite{Kr2},
Lemma~3.9 allow to conclude.
$\Box$

\smallskip

Now consider the Chinese character $S_k$ with $k$ univalent vertices
from the left-hand side of Figure~\ref{SunDiag}.
By definition the element $\chib(S_k)$ is the sum over the $k!$ diagrams
in $\A$ that arise when the univalent vertices of $S_k$ are glued into an
oriented circle in a permuted way.
It follows from the (STU)-relation in $\A$, as
shown in Figure~\ref{STUrel}, that
the element $\chib(S_k)$ can be written as

\begin{equation}\label{SkTk}
\mbox{$\chib(S_k)=k! T_k$ + diagrams that have
$k-1$ vertices on the oriented circle.}
\end{equation}

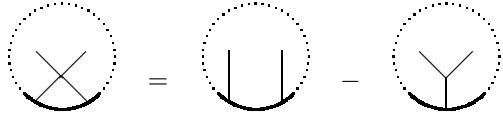
\begin{figure}[!h]
$$
\begin{picture}(2,1.6)(-1,-0.6)
\DottedCircle
\Arc[8]
\Arc[9]
\Arc[10]
\Endpoint[8]
\Endpoint[10]
\thinlines
\put(0.5,-0.866){\line(-1,1){0.966}}
\put(-0.5,-0.866){\line(1,1){0.966}}
\end{picture}\quad = \quad
\begin{picture}(2,1.6)(-1,-0.6)
\DottedCircle
\Arc[8]
\Arc[9]
\Arc[10]
\Endpoint[8]
\Endpoint[10]
\thinlines
\put(0.5,-0.866){\line(0,1){0.966}}
\put(-0.5,-0.866){\line(0,1){0.966}}
\end{picture} \quad - \quad
\begin{picture}(2,1.6)(-1,-0.6)
\DottedCircle
\Arc[8]
\Arc[9]
\Arc[10]
\Endpoint[9]
\thinlines
\put(0.0,-1.0){\line(0,1){0.6}}
\put(0.0,-0.4){\line(1,1){0.5}}
\put(0.0,-0.4){\line(-1,1){0.5}}
\end{picture}
$$
\caption{The STU-relation}\label{STUrel}
\end{figure}

We will now state an application of Theorem~\ref{highestcoeff}
that will be needed later.

\begin{prop}\label{calcd21}
Let $k\geq 4$ be even and let $S_k$ be the Chinese character of degree~$k$
shown on the left-hand side of Figure~\ref{SunDiag}. Then there exist
finite-dimensional
irreducible representations~$\rho_\alpha$ of the Lie superalgebras~$\D21$ such
that for all but finitely many values of~$\alpha$ we have

$$
W_{\D21,\rho_\alpha,\omega_\alpha}(\chib(S_k))\not=0.
$$
\end{prop}

{\bf Proof:}
In Formula~(\ref{D21deltaplus}) we have chosen a positive system
$\Delta^+$ of roots
such that the generators $e_i$ of $\D21$ belong to positive root spaces.
By Theorem~8 of \cite{Kac} there exists a
sequence of finite-dimensional
irreducible representations $\rho_{\alpha,n}$ of $\D21$ with highest
weights $\lambda_{n}=n(3H_1^*+H_2^*+H_3^*)$ ($n=1,2,\ldots$).
By Lemma~\ref{poly2} the values
$W_{\D21,\rho_{\alpha,n},\omega_\alpha}(\chib(S_k))$ form a polynomial
in~$n$ and~$\alpha$ and the degree in~$n$ is bounded by~$k$. Let us
denote the coefficient of $n^k$ by $d_{\alpha,k}$.
By Lemma~\ref{poly2} and Formula~(\ref{SkTk}) we have
$d_{\alpha,k}(\chib(S_k))=k! d_{\alpha,k}(T_k)$.

The bilinear form $<,>_\alpha$ induces a bilinear form on
$H^*$ with matrix

\begin{equation}
\diag((1+\alpha)/2,-1/2,-\alpha/2).
\end{equation}

Now we use Theorem~\ref{highestcoeff} to
compute for even $k\geq 2$ and $\alpha=1$ that

\begin{equation}
d_{1,k}(T_k) =
2\sum_{\beta\in\Delta^+}
(-1)^{\deg\,\beta}<\lambda_{1},\beta>_1^k
= 2 (6^k+2-4^k-2(3^k)-2^k).
\end{equation}

Since we can show that $d_{1,k}(T_k)>0$ for $k\geq 4$,
the coefficient $k! d_{\alpha,k}(T_k)$ of~$n^k$ in
$W_{\D21,\rho_{\alpha,n},\omega_\alpha}(\chib(S_k))$ does not vanish as
a polynomial in~$\alpha$. Thus we may choose~$n_0$ such that
$W_{D_\alpha,\rho_{\alpha,n_0},\omega_\alpha}(\chib(S_k))$ does not vanish
as a polynomial in~$\alpha$.
So with the choice $\rho_\alpha:=\rho_{\alpha,n_0}$
we only have to exclude finitely many values of~$\alpha$ in the statement of
the proposition.
$\Box$

\section{The Elements $D_{k+d,k}$}

Vogel defined a commutative graded algebra~$\Lambda$ equipping the
space~$P(\A)_{\geq 2}$ of pri\-mi\-tive elements of degree~$\geq 2$ with the
structure of a graded $\Lambda$-module.
Let~$\L$ be a simple Lie superalgebra with Casimir~$\omega$.
Then by Theorem~6.1 of~\cite{Vog} there
exists a homomorphism $\chi_{L,\omega}^{}:\Lambda\longrightarrow \C$
such that for all $\lambda\in\Lambda$, for all representations~$\rho$
of~$L$ and for all $a\in P(\A)_{\geq 2}$ we have

\begin{equation}\label{chiL}
W_{\L,\omega,\rho}(\lambda a)=
\chi^{}_{\L,\omega}(\lambda)W_{\L,\omega,\rho}(a).
\end{equation}

In the case of $\D21$ the specialization
$(a,b,c)\mapsto (-\alpha-1,1,\alpha)$ maps
the elements $\sigma_2=ab+ac+bc$ and $\sigma_3=abc$ to
$-1-\alpha-\alpha^2$ and $-\alpha-\alpha^2$ respectively.
By~\cite{Vog}, Theorem~6.13 there exists a morphism of graded algebras
$\chi^{}_D:\Lambda\longrightarrow \C[\sigma_2,\sigma_3]$ $(\deg\,\sigma_i=i)$
such that

\begin{equation}\label{chiD}
\chi^{}_{\D21,\omega_\alpha}(.)
=\chi^{}_D(.)(-1-\alpha-\alpha^2, -\alpha-\alpha^2).
\end{equation}

Certain elements $t,x_3,x_5,x_7,\ldots$ generate a subalgebra
of $\Lambda$
on which formulas for the maps $\chi^{}_{L,\omega}$ are known. These
formulas allow to prove the following proposition.

\begin{prop}\label{calc}
%
For $d=15$ and all $d\geq 17$ there exist elements
$P_d\in \Lambda$ with $\deg\, P_d=d$
that satisfy $\chi^{}_{\L,\omega}(P_d)=0$ for all simple Lie
algebras~$L$
with Casimir~$\omega$, but $\chi^{}_D(P_d)\not=0$.
\end{prop}

{\bf Proof:}
Let $\varphi:\C[T, X_3, X_5, \ldots]\longrightarrow \Lambda$ be the
algebra morphism
defined by mapping~$X_i$ to~$x_i$
and~$T$ to~$t$. Let~$\C[\lambda,\mu,\nu]^{S_3}$ be the ring of
symmetric polynomials in the indeterminates $\lambda,\mu,\nu$.
There exist morphisms of algebras $\chi^{}_0,\chi'_{L,\omega}, \chi'_D$
that make the following diagram commutative
(see~\cite{Vog}, proof of Lemma~7.5).

$$
\begin{picture}(11,4.5)(-5.5,-2)
\put(-5.5,1){\makebox(3,1){$\C[T,X_3,X_5,\ldots]$}}
\put(3,1){\makebox(2,1){$\C[\lambda,\mu,\nu]^{S_3}$}}
\put(-4.5,-2){\makebox(1,1){$\Lambda$}}
\put(3.5,-2){\makebox(1,1){$\C$}}
\put(-1,-0.5){\makebox(2,1){$\C[\sigma_2,\sigma_3]$}}
\put(-2,1.5){\vector(1,0){4.5}}
\put(-4,1){\vector(0,-1){2}}
\put(-3.5,-1.5){\vector(1,0){7}}
\put(4,1){\vector(0,-1){2}}
\put(-3.5,-1){\vector(3,1){2.2}}
\put(3.5,1){\vector(-3,-1){2.2}}
\put(-4.6,-0.5){\makebox(0.6,1){$\varphi$}}
\put(4,-0.5){\makebox(1.5,1){$\chi'_{L,\omega}$}}
\put(-0.75,-1.5){\makebox(1.5,0.7){$\chi^{}_{L,\omega}$}}
\put(-0.5,1.5){\makebox(1,0.7){$\chi^{}_0$}}
\put(-3,-0.8){\makebox(1.2,1){$\chi^{}_D$}}
\put(1.8,0.5){\makebox(1.2,1){$\chi'_D$}}
\end{picture}
$$

It is not difficult to see that the
image of $\chi^{}_0$ is

$$
\C[t]\oplus(t+\lambda)(t+\mu)(t+\nu)\C[t,\lambda\mu+\mu\nu+\nu\lambda,
\lambda\mu\nu],
$$

where $t=\lambda+\mu+\nu$.
Let $I=\{\lambda,\mu,\nu\}$.
Let

$$ P= \prod_{a\in I}(t+a)
\prod_{a\in I}(t-a)
\prod_{a,b\in I \atop a\not=b}(a+2b) 
\prod_{a\in I}(3a-2t)
\in\C[\lambda,\mu,\nu]^{S_3}{}. $$

For all $Q\in\C[\lambda,\mu,\nu]^{S_3}$ the element~$PQ$
is in the image of~$\chi^{}_0$. Furthermore,
we have $\chi'_{L,\omega}(PQ)=0$ for all simple Lie algebras~$L$ with
Casimir~$\omega$ (see the computation of eigenvalues of~$\Psi$ on
$S^2 L/\omega$ in \cite{Vog}).

Let~$Q\in\C[\lambda,\mu,\nu]^{S_3}$ be a homogeneous element not
divisible by~$t$.
Then we have $\chi'_D(PQ)\not=0$.
We have~$\deg(PQ)=15+\deg\, Q$ and $\deg\, Q=0$ or $\deg\, Q\geq 2$.
Define~$P_d=\varphi(p)$, where
$p\in \chi_0^{-1}(PQ)$ and $d=\deg\, p=\deg(PQ)$.
The element~$P_d$ has the properties stated in the proposition.$\Box$

\medskip

Proposition~\ref{calc} was proved by the author by making a
computation with Mathematica (see~\cite{Li2}). The proof stated above
is due to P.\ Vogel.

\medskip

{\bf Proof of Theorem~\ref{ineq}:}
We turn the space $P(\B)_{\geq 2}$ of primitive elements of $\B$
of degree $\geq 2$
into a $\Lambda$-module by the formula

\begin{equation}
\chib(\lambda \cdot b)=\lambda\cdot \chib(b) \quad
(\lambda\in\Lambda, b\in\B).
\end{equation}

An example is shown in Figure~\ref{bigdiag}.


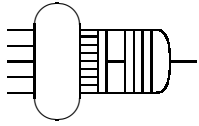
\begin{figure}[!h]
{\unitlength=17pt
$$
\begin{picture}(5,2.8)(-2.5,-0.8) 
\thinlines
\put(-0.4,-0.7){\line(1,0){1.7}}
\put(-0.4,0.7){\line(1,0){1.7}}
\put(2.2,0){\line(-1,0){0.6}}
\qbezier(1.3,0.7)(1.6,0.7)(1.6,0.0)
\qbezier(1.3,-0.7)(1.6,-0.7)(1.6,0.0)
\multiput(1.2,-0.7)(-0.2,0){3}{\line(0,1){1.4}}
\multiput(0.6,-0.7)(-0.4,0){2}{\line(0,1){1.4}}
\put(0.6,0.0){\line(-1,0){0.4}}
\multiput(0,-0.7)(-0.4,0){2}{\line(0,1){1.4}}
\multiput(0,-0.53)(0,0.18){7}{\line(-1,0){0.4}}
\put(-0.9,0){\oval(1,2.6)}
\multiput(-2,-0.7)(0,0.35){5}{\line(1,0){0.6}}
\end{picture}
$$
\caption{The element $t^3 x_3 x_9\cdot S_6$.}\label{bigdiag}
}
\end{figure}

Now we can define $D_{k+d,k}=P_d\cdot S_k$
with the elements $P_d$ of Proposition~\ref{calc}.
By Formula~(\ref{chiL}) we have
$W_{\L,\omega,\rho}(\chib(D_{k+d,k}))=0$ for all simple Lie algebras~$L$ with
Casimir~$\omega$ and representation~$\rho$.
By standard arguments using
Exercise~6.33 of~\cite{BN1}
this can be generalized to the case where~$L$ is a semisimple Lie algebra.
Let $\rho_\alpha$ be the
representation of~$\D21$ from Proposition~\ref{calcd21}. Then by
Formula~(\ref{chiL}) and Formula~(\ref{chiD}) we have

\begin{equation}
\WDroa(\chib(D_{k+d,k}))=\chi^{}_D(P_d)
(-1-\alpha-\alpha^2,-\alpha-\alpha^2)\WDroa(\chib(S_k)).
\end{equation}

By Proposition~\ref{calcd21}
the factor $\WDroa(\chib(S_k))$ only vanishes for finitely many
values of $\alpha$.
With our choice of $P_d$ the polynomial $\chi^{}_D(P_d)
\in \C[\sigma_2,\sigma_3]$ does not vanish and is homogeneous with respect to
$\deg\, \sigma_i=i$. So the value
$\chi^{}_D(P_d)(-1-\alpha-\alpha^2,-\alpha-\alpha^2)$ can also only vanish
for finitely many choices of $\alpha$. This implies $D_{k+d,k}\not=0$ and
completes the proof.
$\Box$

\section*{Acknowledgements}
I thank C.~Kassel and P.~Vogel for carefully reading
the paper and for proposing improvements.
Thanks to the German Academic Exchange Service for financial support
(Doktorandenstipendium HSP~III).
}


\begin{thebibliography}{4}
\bibitem
{BN1}
D.\ Bar-Natan, {\em On the Vassiliev knot invariants}, Topology 34
(1995), 423--472.

\bibitem
{BN2} D.\ Bar-Natan, {\em Some Computations Related to Vassiliev
Invariants}, available at http://www.ma.huji.ac.il/\~{}drorbn, May 5, 1996.

\bibitem
{BNG} D.\ Bar-Natan and S.\ Garoufalidis, {\em On the Melvin-Morton-Rozansky
conjecture},
Invent.\ Math.\ 125 (1996), 103--133.

\bibitem
{Hum} J.\ E.\ Humphreys, {\em Introduction to Lie Algebras and
Representation Theory}, Springer-Verlag, 1989.

\bibitem
{Kac} V.\ C.\ Kac, {\em A sketch of Lie superalgebra theory},
Comm.\ Math.\ Phys. 53 (1977), 31--64.

\bibitem
{Kne} J.\ A.\ Kneissler, {\em The number of primitive
Vassiliev invariants up to degree twelve}, q-alg/9706022 and
University of Bonn preprint, June 1997.

\bibitem
{Kr2} A.\ Kricker, {\em Alexander-Conway limits of many
Vassiliev weight systems},
J.\ of Knot Th.\ and its Ramif., Vol.\ 6, Nr.\ 5 (1997), 687--714.

\bibitem
{Kr1} A.\ Kricker, B.\ Spence and I.\ Aitchison,
{\em Cabling the Vassiliev invariants},
J.\ of Knot Th.\ and its Ramif., Vol.\ 6, Nr.\ 3 (1997), 327--358.

\bibitem
{Li2} J.\ Lieberum, {\em charcomp.m}, a Mathematica program to
perform computations with the characters from \cite{Vog}, available at
http://www-irma.u-strasbg.fr/irma/, 1997.

\bibitem
{Li1} J.\ Lieberum, {\em Chromatic weight systems and the
corresponding knot invariants}, preprint I.R.M.A.
Strasbourg (1997), to appear in Math.\ Ann.

\bibitem
{Vog} P. Vogel, {\em Algebraic structures on modules of
diagrams}, Universit\'{e}
Paris VII preprint (1995), to appear in Invent.\ Math.
\end{thebibliography}
\end{document}